\documentclass[authoryear,preprint,11pt]{elsarticle}
\usepackage{amssymb}
\usepackage{eurosym}
\usepackage{color,colortbl}
\usepackage{bm}
\usepackage{booktabs}
\usepackage{multirow}
\usepackage{longtable}
\usepackage{arydshln}
\usepackage{comment}
\usepackage{tabu}
\usepackage[a4paper, total={8in, 8in}]{geometry}
\usepackage{caption} % To customize the captions
\usepackage{subcaption}
\usepackage{tikz}
\usepackage{pgfplots}
\usetikzlibrary{arrows,calc}
\usepackage{graphicx}
\usepackage{fancybox}
\usepackage{pifont}
\usepackage{multirow}
\usepackage{lscape}
\usepackage{enumitem}

\usepackage{amsthm,amsmath,amssymb,eurosym}

\makeatletter
\def\ps@pprintTitle{%
  \let\@oddhead\@empty
  \let\@evenhead\@empty
  \def\@oddfoot{\reset@font\hfil\thepage\hfil}
  \let\@evenfoot\@oddfoot
}
\makeatother

\journal{TBC}

\usepackage{lineno, blindtext}

\begin{document}

\begin{frontmatter}
    \title{A Multi-Objective Portfolio of Portfolios Problem with \\ Qualitative Performance Assessments}
    \author[corl]{{\sc M. Barbati}\corref{cor1}}\ead{maria.barbati@unive.it}
    \author[deb]{{\sc S. Greco}}
    \author[ist]{{\sc J.R. Figueira}}
   \address[corl]{Deparment of Economics, Ca' Foscari University of Venice, Italy}
    \address[deb]{Department of Economics and Business, University of Catania, Italy}
    \address[ist]{CEGIST, Instituto Superior T\'{e}cnico,  Universidade de Lisboa, Portugal}

    \cortext[cor1]{Corresponding author at: Deparment of Economics, Ca' Foscari University of Venice, Italy}

    \begin{abstract}
   \noindent {We present a multi-objective portfolio decision model that involves selecting both a portfolio of projects and a set of elements to allocate to each project. Our model includes a defined set of objectives to optimize, with projects contributing to these objectives in various ways. The elements included in the portfolios are assessed based on both qualitative and quantitative criteria. Projects can only be selected for the portfolio if they meet specific requirements defined by threshold values on the criteria.  The model is adaptable to include temporal considerations and stochastic, making it suitable for a wide range of real-life applications. To manage the decision-making process, we employ an interactive multi-objective method that integrates the selection of both portfolios and elements. After making initial selections, we ask the decision-maker to evaluate the portfolios, from which we derive a series of rules to be incorporated into the multiobjective model until the decision-maker is satisfied. We illustrate the functionality of our model through an illustrative case study.}
    \end{abstract}
    \vspace{0.25cm}
\end{frontmatter}

%%\vfill\newpage

%%\tableofcontents

%%\vfill\newpage

    %%
    \vspace{0.25cm}
\section{Introduction}\label{sec:introduction}
\noindent {Multi-objective portfolio models have been developed to select a set of projects based on multiple objectives \citep{SaloEtAl11}. To implement these projects, various elements are required, such as human resources or capital. However, only a limited number of these elements can be selected, and they often need to be shared among multiple projects. Furthermore, elements can be evaluated according to specific criteria, meaning that only those elements meeting a certain level of ``quality" can be assigned to individual projects. As a result, the selection of the portfolio of projects can depend on both the choice of elements and their ``quality."}

{In this paper, we propose a model that examines the relationship between portfolio selection and the selection of individual elements within that portfolio. Specifically, our objectives are to achieve the following:
\begin{itemize}[label={--}]
\item Develop a portfolio of projects that optimizes certain objectives.
\item Create a portfolio of elements to be assigned to each selected project, ensuring that these elements meet specified thresholds for various criteria.
\end{itemize} 
}

{The ``quality" of the projects within the portfolio is determined by the ``quality" of the portfolios of elements allocated to each project. Specifically, a project can only be selected to contribute to the optimization objectives if enough elements are assigned that exceed a specified threshold for each criterion.}

{Let us explore some potential applications of our model. As a first example, consider a company that aims to select a portfolio of Research and Development (R\&D) projects. The company must decide which projects to invest based on strategic objectives, such as the scientific results to be achieved. In addition, researchers need to be assigned, through specific portfolios, to each project. 
These portfolios of researchers must meet defined thresholds for certain criteria. For example, to run a specific project, the company requires at least two researchers, each with at least three years of experience in R\&D projects. Consequently, the company must determine which portfolio of projects to select and which portfolios of researchers to be assigned to each project. In essence, they need to choose a portfolio composed of R\&D projects that include portfolios of researchers, ensuring that all requirements are satisfied and optimizing various objectives.}

{As a second example, consider a urban planning problem. Suppose a council needs to choose a portfolio of services to offer to citizens while optimizing objectives such as quality and accessibility. They require a set of organizations (our elements) to provide these services. These organizations are selected based on specific criteria such as profitability, reliability, and environmental impact. For example, the council might stipulate that the organizations responsible for waste collection cannot exceed a predetermined threshold for environmental impact.}

{As a third example, how can a retailer select a portfolio of stores to store? With increasing budget constraints, many retailers face the challenge of determining which stores to remain open based on objectives such as sales performance and the attractiveness of their locations. Simultaneously, they must decide which personnel to retain, evaluating their skills against the requirements needed to manage each shop effectively.}

{It is important to note that the examples mentioned above are just a few potential applications. Additional examples could include selecting a portfolio of suppliers based on the range of services they can provide or scheduling and selecting health services for patients.}

{Although there are examples in the literature that address portfolio problems involving the allocation of elements to projects \citep[see e.g.][]{gutjahr2008competence}, to the best of our knowledge, our model is the first to introduce the selection of a portfolio of elements based on specific threshold levels. This capability allows the model to accommodate both qualitative and quantitative evaluations. Our model features the following innovative characteristics:}

{
\begin{itemize}[label={--}]
\item Portfolios of projects are selected based on strategic objectives and the ``quality" provided by the portfolios of elements allocated to each project, ensuring that they meet certain threshold levels.
\item A portfolio contributes to the overall objectives if specific conditions related to its elements are satisfied. For example, if a certain number of elements meet the required level for a particular criterion, then the project is deemed to contribute to the specified value for that objective.
\item Our model allows for considering various objectives and criteria in selecting both portfolios and elements. It effectively handles both qualitative and quantitative evaluations without requiring the aggregation of evaluations across multiple criteria.
\item Multiple logical restrictions can be defined and incorporated as constraints within our model.
\end{itemize}
}

{
 In addition to the points mentioned previously, we enhance our model by considering the temporal distribution of selected projects and the uncertainty associated with various potential scenarios. In many applications, it is essential to distribute the chosen projects throughout the planning period, ensuring that the available resources meet the defined thresholds at each stage. Furthermore, it is often necessary to incorporate uncertain evaluations for both the strategic objectives related to the project portfolio and the criteria that must be satisfied by the resources allocated to these projects. 
}

{
Several methodologies can be implemented to manage our model effectively. We demonstrate how it can be handled through interactive methodologies adapted from \citep{barbati2018}, which benefit from stakeholder interactions and the learning processes that stakeholders undergo during these interactions \citep{Miettinen2016}.
}

{
The paper is organized as follows. Section \ref{sec:literature}, reviews the literature on multi-objective portfolio selection. Section \ref{sec:mo-portfolios} presents the formulation of our portfolio of portfolios model, along with two possible extensions: the temporal distribution of the projects and the uncertainty of conditions. Section \ref{sec:methodology}, introduces the methodology used to address the problem. Section \ref{sec:bigger_example} explains how to apply an interactive multi-objective approach to manage the model and its possible extensions, while Section \ref{sec:conclusion} concludes the paper.
}

\section{Literature review}\label{sec:literature}
%\section{Concepts, definitions, and notation}\label{sec:concepts}

\noindent {Portfolio decision analysis is a strand of research focused on selecting a set of projects based on specific criteria while satisfying various constraints \citep{salo2011portfolio}. The simplest scenario involves optimizing a single objective, such as the Net Present Value (NPV) of the selected projects \citep{ball1969portfolio}. A more complex and interesting approach considers multiple objectives that may vary depending on the application context \citep{liesio2020portfolio}. In such a case, several methods can be utilized, including exact algorithms \citep{liesio2007preference}, preference methods \citep{zheng2011constrained}, and interactive methods \citep{barbati2018}. Additionally, various factors can be addressed, such as uncertainty \citep{relich2017fuzzy}, interactions among projects \citep{liesio2008robust}, and the involvement of multiple stakeholders \citep{gade2018redis}.}

{When selecting projects, decisions often involve not only which projects to choose but also the timing of these projects \citep{ghasemzadeh2000project}. This area of research, commonly referred to as project selection and scheduling \citep{carazo2010solving}, has seen the development of various models that consider multiple factors \citep{perez2016multiobjective}. These problems have been addressed using a range of methods, including heuristic approaches, as demonstrated by \citep{huang2016uncertain}, as well as exact methods like the $\epsilon$-constraint method \citep{rezahoseini2020comprehensive}.}

{When defining the schedule for projects, a crucial aspect is the allocation of shared resources among the various projects. This area of research is often referred to as resource-constrained project scheduling. Several objectives can be considered in this context, such as minimizing lateness, tardiness, or earliness of the projects \citep{hartmann2010survey}. Typically, these objectives focus on individual projects; however, in some cases, they extend to consider the entire project portfolio \citep{browning2016managing}. Additionally, various complexities arise from uncertainties in evaluating the available resources, which can be modeled using fuzzy attributes \citep{mohanty2005fuzzy}. The situation becomes even more challenging when the shared resources involve human skills \citep{barreto2008staffing}.}

{In particular, \citep{gutjahr2008competence} developed a model that selects a set of projects and allocates available human resources over a specified time horizon, taking into account the evolution of efficiency and the improvement of competencies during that period. This model was further enhanced in a subsequent study where the authors considered a multi-objective version of the problem and proposed several metaheuristic methods to address it \citep{gutjahr2010multi}. Additionally, in another work,  \citep{gutjahr2010bi} introduced uncertainty in the evolution of competencies over time and presented a new methodology to exploit the entire Pareto front. Furthermore, the authors investigated the theoretical properties that could encourage a company to invest in projects of the same class, meaning those with similar characteristics \citep{gutjahr2011optimal}. In a later study \citep{gutjahr2013project},  they described outsourcing costs and introduced a new two-step approach.}

{In the same strand of research,  \citep{zaraket2014skill}  developed a framework to classify the software skills possessed by personnel and introduced a portfolio selection model that allocates resources from various classes of software skills. A  similar issue was addressed by \citep{hematian2020multiobjective}, who considered the duration of activities associated with the learning process. 
Various methodologies have been employed to tackle these problems, including a nonlinear mixed-integer programming approach \citep{gutjahr2008competence}  and metaheuristic methods, such as the combination of NSGA-II with a Pareto ant colony optimization algorithm \citep{chen2020competence}. However, to the best of our knowledge, only \citep{lukovac2017portfolio} incorporated the preferences of decision-makers in an interactive decision process. Their work created all potential assignment scenarios for human resources based on a defined set of constraints.}

{It is important to note that the approach taken by this latter strand of research differs from ours. While they do select a portfolio of projects and assign human resources to them, they do not determine which human resources should be retained or removed. In other words, their assumption is that the available resources are predetermined, and their decision-making process only concerns how these resources are allocated among the various projects.}

{
In a recent development in this field \citep{noro2023project}, the authors focus solely on selecting projects for portfolios, treating agents (i.e., the elements allocated to each portfolio) as fixed and only distributing them, rather than selecting them. Our paper aims to address this gap by examining both project and elemets selection.}

\section{The multi-objective portfolio of portfolios problem}\label{sec:mo-portfolios}

\noindent In this section, we present the problem of a portfolio of portfolios, a scenario commonly encountered in large investment companies where multiple individual portfolios are combined into a larger, overarching portfolio. This primary problem is examined alongside two significant extensions: incorporating temporal considerations and integrating stochastic elements.  For each of these three problems, we provide illustrative example to clarify the concepts  thereby offering a better understanding of the three problems.

\subsection{Problem statement}\label{subsec:pb_mo-portfolios}
\noindent Let us consider a set of projects, $P = \{P_1, …, P_j,…,P_n\}$, where each project leads to a certain benefit, say the net present value, and needs a given amount of investment to be implemented. Each project can be selected or not. Knowing that a budgetary constraint exists, it is impossible to select all the projects in $P$ simultaneously. The classical portfolio selection problem consists of selecting a subset of $P$, called the \textit{portfolio}, which maximizes the overall net present value and fulfills the budgetary constraint. In practice, in general, we are faced with other benefits than the net present value, such as, for example, the number of workers hired if a project is selected and the environmental impact of the project. In the latter case, we are in the presence of a \textit{multi-objective portfolio selection problem}. 

Our problem, called \textit{portfolio of portfolios (PoP)}, is an extension of the previous one with a two-fold characteristics, one is that each project $P_j$, for $j=1,\ldots,n$, is formed from a set of components or elements, $E=\{e_1,\ldots,e_i,\ldots,e_m\}$, the other is that the performances of each element, $e_i$, for $i=1,\ldots,m$, are assessed according to several criteria, $g_p$, for $p=1,\ldots,q$, of a quantitative and/or qualitative nature. In addition, we will consider the multi-objective $PoP$  problem.

The following additional notation is needed for the rest of the paper. Let,

\begin{itemize}[label={--}]
    \item $c^{\ell}_j$ denote the value of project, $P_j$, according to the objective function, $z_\ell$, for $j=1,\ldots, n$, and $\ell = 1,\ldots, k$.  
    \item $g_p(e_i)$ denote the performance level of element $e_i$ according to criterion $g_p$, for $i=1,\ldots,m$ and $p=1,\ldots q$.
    \item $L_p = \{l_{p1},\ldots, l_{ph}, \ldots, l_{pf{_p}}\}$ denotes the quantitative or qualitative discrete scale of criterion $g_p$, for $p=1,\ldots,q$.
    \item $v_{iph}$ denotes a parameter with a value equal to $1$ if $g_p(e_i) \geqslant l_{ph}$, and $0$ otherwise, for $i=1,\ldots,m$, $p=1,\ldots,q$, and $h=1,\ldots,f_p$. 
    \item $u_{jph}$ denotes the minimal number of elements, $e_i \in E$, assigned to project, $P_j \in P$, for which it is required that $g_p(e_i) \geqslant l_{ph}$; $D$ is the set containing all triplets $(j,p,h)$. 
    \item $w_{ij}$ denoted the investment or the cost of component of element $e_i$, when assigned to project, $P_j$, for $i=1,\ldots, n$, and $j=1,\ldots,n$.
    \item $W$ denotes the available budget for the investment in all the assigned elements to form the selected projects. 
\end{itemize}

Two types of decision variables can be defined as follows:

\begin{itemize}[label={--}]
    \item $x_j$, which is equal to $1$ if project, $P_j$, is selected, and $0$ otherwise, for $j=1,\ldots,n$. 
    \item $y_{ij}$, which is equal to $1$ if component, $e_i$, is assigned to project, $P_j$, and $0$, otherwise, for $i=1,\ldots,m$ and $j=1,\ldots,n$.
\end{itemize}

\begin{equation}\tag{\textit{MO-PoP}}
    \begin{array}{lcl}
         {\displaystyle \max \; z_1(x_1,\ldots,x_j,\ldots,x_n) } & = & {\displaystyle \sum_{j=1}^{n}c^{1}_jx_j}, \\
         & \vdots & \\
        {\displaystyle \max \; z_\ell(x_1,\ldots,x_j,\ldots,x_n) } & = & {\displaystyle \sum_{j=1}^{n}c^{\ell}_jx_j}, \\
         & \vdots & \\
          {\displaystyle \max \; z_k(x_1,\ldots,x_j,\ldots,x_n) } & = & {\displaystyle \sum_{j=1}^{n}c^{k}_jx_j},\\
         \mbox{subject to:} & & \\
         & & {\displaystyle \sum_{j=1}^{n}y_{ij} \leqslant 1, \;\;\; i = 1, \ldots, m,} \\
         & & {\displaystyle y_{ij} \leqslant x_j, \;\;\; i = 1, \ldots, m, \; j = 1, \ldots, n,} \\
         & & {\displaystyle \sum_{i=1}^{m}v_{iph}y_{ij} \geqslant u_{jph}x_j, \;\;\; (j,p,h) \in D,} \\
         & & {\displaystyle \sum_{i=1}^{m}\sum_{j=1}^{n}w_{ij}y_{ij} \leqslant W.} \\
    \end{array}
\end{equation}

The first set of constraints models that an element, $e_i$, will be assigned at most to one project $P_j$. If an element, $e_i$, is by default an element of a certain project, $P_j$, the variable $y_{ij}$ will be \textit{a priori} set to one. The second set of constraints models that if a project is not selected, it does not contain any elements or components. The third set of constraints models the minimal required number of elements assigned to a given project with a certain minimal performance level requirements. The last constraint is the budgetary constraint.

It is easy to see this model can be extended to one  with more levels, if each components can itself be composed of sub-components and so on. 

%\textcolor{black}{It is important to note that while we are treating the allocation of a project as a binary variable, we could also consider including only specific components of a project. In this case, we would use continuous variables to represent the portions of the projects being included in the portfolio.}

\subsection{An illustrative example}\label{subsec:ex_mo-portfolios}
\noindent Let us consider a company aiming to determine the most suitable R\&D portfolio of projects to implement in the next planning period. The company aims to maximize the following three strategy objectives: the scientific quality of a portfolio of projects ($z_1$); the adequacy of a portfolio to the organization strategy ($z_2$); and, the impact generated by a portfolio outside the organization, such as the innovation developed and the benefits for the population ($z_3$). 

Each project is composed of researchers that must be assigned to it. In accordance with the internal regulations, the company has defined three gain-type criteria to assess the performance of the researchers, as follows: the number of months of experience in managing R\&D projects of each researcher ($g_1$); a qualitative score on a 3-point scale for assessing the expertise of each researcher ($g_2$); and, the number of tasks each researcher completes on time ($g_3$).

Assume we have four researchers, i.e., $E=\{e_1,e_2,e_3,e_4\}$, and three projects, i.e., $P=\{P_1,P_2,P_3\}$. Table 1.a presents the performance levels of each researcher on each criterion. For the sake of simplicity, we define only two levels for each scale: $L_1 = \{l_{11}=20, l_{12}=40\}$, $L_2 = \{l_{21}=2, l_{22}=3\}$, and $L_3 = \{l_{31}=20, l_{32}=40\}$. Now, we can compute the values of the parameter, $v_{iph}$, and presented in Table 1.b. For example, $v_{111}=0$ since $g_1(e_1) = 18$ is strictly lower than a given $l_1 = 20$. Instead, $v_{211} = 1$ since $g_1(e_2)=24$ is greater than of equal to $l_1 = 20$.   

\begin{table}[htpb!]
  \begin{minipage}{.3\linewidth}
    \centering
    \begin{tabular}{cccc} \toprule 
            & $g_1$ & $g_2$ & $g_3$ \\ \hline
      $e_1$ &   18  &   3   & 42    \\
      $e_2$ &   24  &   2   & 38    \\ 
      $e_3$ &   43  &   1   & 54    \\ 
      $e_4$ &   36  &   3   & 49    \\ \bottomrule
    \end{tabular}
    \caption*{a. Performance levels $g_p(e_i)$.}
  \end{minipage}%
  \hspace{1.25cm}
  \begin{minipage}{.3\linewidth}
    \centering
    \begin{tabular}{ccccccccc}
    \toprule
    \multicolumn{1}{c}{} & \multicolumn{2}{c}{$g_1$} & \multicolumn{2}{c}{$g_2$} & \multicolumn{2}{c}{$g_3$}  \\
    \cmidrule(rl){2-3} \cmidrule(rl){4-5} \cmidrule(rl){6-7} 
    $v_{iph}$ & $l_{11}$ & $l_{12}$ & $l_{21}$ & $l_{22}$ & $l_{31}$ & $l_{32}$  \\
    \midrule
    $e_1$ & 0 & 0 & 1 & 1& 1 & 1  \\
    $e_2$ & 1 & 0& 1 & 0 & 1 & 0 \\
    $e_3$ & 1 & 1 & 0 & 0 & 1 & 1  \\
    $e_4$ & 1 & 0 & 1 & 1 & 1 & 1  \\
    \bottomrule
    \end{tabular}
    \caption*{b. Parameters $v_{iph}$.}
  \end{minipage}%
  \caption{Performance levels of the elements and minimal qualitative levels required.}
\end{table}

Each researcher will earn a salary that may be different according to the project she/he will be assigned to. The salary is in K\euro{} and provided in Table 2.a. The available budget is 100 K\euro{}. 

The values of parameter $u_{jph}$ are displayed in Table 2.b. For example, the value $u_{111} = 2$ means that at least two elements that have performances levels on $g_1$ at least equal to $l_{11} = 20$  should be selected for project $P_1$, i.e., at least two researchers that have at least 20 months of experience in R\&D development will be necessary for implementing project $P_1$. Note that all these values are not defined for all the triplets $(j,p,h)$, for example, $u_{112=0}$, which means that there is no constraint on the number of elements to be selected for project $P_1$. We do not need to include any element in project $P_1$ that has a performance level for criterion $g_1$ at least equal to $l_{12}=40$, i.e., to implement project $P_1$ we do not need any researcher with more than $40$ months of experience.

%%\vfill\newpage

\begin{table}[htpb!]
  \begin{minipage}{.3\linewidth}
    \centering
    \begin{tabular}{cccc} \toprule
            & $w_{i1}$ & $w_{i2}$ & $w_{i2}$ \\ \hline
      $e_1$ &   37  &   35  & 30    \\
      $e_2$ &   40  &   35  & 33    \\ 
      $e_4$ &   44  &   26  & 36    \\ 
      $e_5$ &   32  &   28  & 34    \\ \bottomrule
    \end{tabular}
    \caption*{a. Research salaries, $w_{ij}$.}
  \end{minipage}%
  \hspace{1.75cm}
  \begin{minipage}{.3\linewidth}
    \centering
    \begin{tabular}{ccccccc}
    \toprule
    \multicolumn{1}{c}{} & \multicolumn{2}{c}{$P_1$} & \multicolumn{2}{c}{$P_2$} & \multicolumn{2}{c}{$P_3$} \\
    \cmidrule(rl){2-3} \cmidrule(rl){4-5} \cmidrule(rl){6-7} 
    $u_{jph}$ & $l_{11}$ & $l_{12}$ & $l_{21}$ & $l_{22}$ & $l_{31}$ & $l_{32}$ \\
    \midrule
    $g_1$ & 2 & 0 & 1 & 0 & 1 & 0  \\
    $g_2$ & 1 & 0 & 1 & 1 & 2 & 0  \\
    $g_3$ & 1 & 0 & 1 & 0 & 1 & 1  \\
    \bottomrule
    \end{tabular}
    \caption*{b. Minimal number of elements assigned to each project, $u_{jph}$.}
  \end{minipage}%
  \caption{Investment costs and minimal number of elements for each project.}
\end{table}

\noindent Finally, Table 3 displays the values of each project according to each objective function. 

\begin{table}[htpb!]
    \centering
     \begin{tabular}{cccc} \toprule
            & $c_{j}^1$ & $c_{j}^2$ & $c_{j}^3$ \\ \hline
      $P_1$ &   24  &   145 & 42    \\
      $P_2$ &   43  &   54  & 24    \\ 
      $P_3$ &   75  &   150 & 57    \\ \bottomrule
    \end{tabular}
    \caption{Objective function values}
\end{table}

\textcolor{black}{We generated a set of solutions on the Pareto Front using the $\epsilon$-constraint method, specifically following the approach proposed by \cite{mesquita2023new}. This algorithm is designed to represent the Pareto front for multi-objective linear programming problems involving two or more objective functions. As a result, we successfully identified nine efficient portfolios along with their corresponding outcome vectors.}

\begin{enumerate}
   
    \item The portfolio of projects $P_1$ and $P_2$ is selected, with correspondent variables $x_1=x_2=1$ and $x_3=0$,
    with the following possible combinations of researchers:
    \begin{itemize}[label={--}]
        \item researchers $e_2$ and $e_3$ allocated to project $P_1$ and researcher $e_4$ allocated to project $P_2$, with variables $y_{21},y_{31},y_{42}$ equals to 1 and all others equal to $0$;
        \item researchers $e_2$ and $e_4$ allocated to project $P_1$ and researchers $e_1$ and $e_3$ allocated to project $P_2$, 
 with variables $y_{21},y_{41},y_{12},y_{32}$ equals to 1 and all others equal to $0$;
        \item researchers $e_
        3$ and $e_4$ allocated to project $P_1$ and researchers $e_1$ and $e_2$ allocated to project $P_2$,  with variables $y_{31},y_{41},y_{12},y_{22}$ equals to 1 and all others equal to $0$;
        \item researchers $e_2$ and $e_3$ allocated to project $P_1$ and researchers $e_1$ and $e_4$ allocated to project $P_2$, with variables $y_{21},y_{31},y_{12},y_{42}$ equals to 1 and all others equal to $0$;
        \item researchers $e_1$, $e_2$ and $e_3$ allocated to project $P_1$ and researcher $e_4$ allocated to project $P_2$, with variables $y_{11},y_{21},y_{31},y_{42}$ equals to 1 and all others equal to $0$;
    \end{itemize}
    \item The portfolios of projects $P_2$ and $P_3$, with variables $x_2=x_3=1$ and $x_1=0$ is selected
    with the following possible combinations of researchers:
    \begin{itemize}[label={--}]
        \item researchers $e_3$ and $e_4$  allocated to project $P_2$ and researcher $e_1$ and $e_2$ allocated to project $P_3$, with variables $y_{32},y_{42},y_{13},y_{23}$ equal to 1 and all others equal to $0$;
        \item researchers $e_1$ and $e_3$  allocated to project $P_2$ and researcher $e_2$ and $e_4$ allocated to project $P_3$, with variables $y_{12},y_{32},y_{23},y_{43}$ equal to 1 and all others equal to $0$;
        \item researcher $e_4$ allocated to project $P_2$ and researcher $e_1$, $e_2$ and $e_3$ allocated to project $P_3$, with variables $y_{42},y_{13},y_{23},y_{33}$ equal to 1 and all others equal to $0$;
        \item researchers $e_4$ allocated to project $P_2$ and researcher $e_2$ and $e_3$ allocated to project $P_3$, with variables $y_{42},y_{23},y_{33}$ equal to 1 and all others equal to $0$.
    \end{itemize}
\end{enumerate}

\subsection{Adding temporal considerations}\label{subsec:temporal_mo-portfolios}

\noindent Let us consider  the same problem under time periods, $t= 1, \ldots, r$. 

\begin{equation}\tag{\textit{tMO-PoP}}
    \begin{array}{lcl}
         {\displaystyle \max \; z_1(x_{t1},\ldots,x_{tj},\ldots,x_{tn}) } & = & {\displaystyle \sum_{t=1}^{r}\sum_{j=1}^{n}c^{1}_{tj}x_{tj}}, \\
         & \vdots & \\
        {\displaystyle \max \; z_\ell(x_{t1},\ldots,x_{tj},\ldots,x_{tn}) } & = & {\displaystyle \sum_{t=1}^{r}\sum_{j=1}^{n}c^{\ell }_{tj}x_{tj}},\\
         & \vdots & \\
        {\displaystyle \max \; z_k(x_{t1},\ldots,x_{tj},\ldots,x_{tn}) } & = & {\displaystyle \sum_{t=1}^{r}\sum_{j=1}^{n}c^{k}_{tj}x_{tj}}, \\
         \mbox{subject to:} & & \\
         & & {\displaystyle \sum_{t=1}^{r}\sum_{j=1}^{n}y_{tij} \leqslant 1, \;\;\; i = 1, \ldots, m,} \\
         & & {\displaystyle y_{tij} \leqslant x_{tj}, \;\;\; t = 1, \ldots, r, \; i = 1, \ldots, m, } \\
         & & {\displaystyle \;\;\;\;\;\;\;\;\;\;\;\;\;\;\;\;\;\;\; j = 1, \ldots, n,} \\
         & & {\displaystyle \sum_{i=1}^{m}v_{tiph}y_{tij} \geqslant u_{tjph}x_{tj}, \;\;\; (t,j,p,h) \in D^t,} \\
         & & {\displaystyle \sum_{t=1}^{r}\sum_{i=1}^{m}\sum_{j=1}^{n}w_{tij}y_{tij} \leqslant W.} \\
    \end{array}
\end{equation}

\subsection{A second illustrative example}\label{subsec:ex_sec_mo-portfolios}

\noindent Let us suppose that the projects described in the previous example can be activated in two different periods, so we have $t=\{1,2\}$. For the sake of the data provided, we suppose that the evaluations $g_p(e_i)$ and the parameters $v_{iph}$ reported in the previous subsection are the evaluations related to the first period, i.e. $g_p^1(e_i)$ and $v_{1iph}$. Instead, in Tables 4.a and 4.b we report the $g_{p}^2(e_i)$ and the parameters $v_{2iph}$ for the second period.

\begin{table}[htpb!]
  \begin{minipage}{.3\linewidth}
    \centering
    \begin{tabular}{cccc} \toprule 
            & $g_1$ & $g_2$ & $g_3$ \\ \hline
     \textit{$e_1$} 	&	42	&	2	&	45		\\
\textit{$e_2$}	&	39	&	2	&	26		\\
\textit{$e_3$}	&	19	&	1	&	22		\\
\textit{$e_4$} 	&	18	&	3	&	35		\\ \bottomrule
    \end{tabular}
    \caption*{a. Performance levels $g_p^2(e_i)$.}
  \end{minipage}%
  \hspace{1.25cm}
  \begin{minipage}{.3\linewidth}
    \centering
    \begin{tabular}{ccccccccc}
    \toprule
    \multicolumn{1}{c}{} & \multicolumn{2}{c}{$g_1$} & \multicolumn{2}{c}{$g_2$} & \multicolumn{2}{c}{$g_3$}  \\
    \cmidrule(rl){2-3} \cmidrule(rl){4-5} \cmidrule(rl){6-7} 
    $v_{iph}$ & $l_{11}$ & $l_{12}$ & $l_{21}$ & $l_{22}$ & $l_{31}$ & $l_{32}$  \\
    \midrule
    $e_1$ & 1 & 1 & 1 & 0 & 1 & 1 \\
$e_2$ & 1 & 0 & 1 & 0 & 0 0 \\
$e_3$ & 0 & 0 & 0 & 0 & 1 & 0  \\  
$e_3$ & 0 & 0 & 1 & 1 & 1 & 0  \\
    \bottomrule
    \end{tabular}
    \caption*{b. Parameters $v_{2iph}$.}
  \end{minipage}%
  \caption{Performance levels of the elements and minimal qualitative levels required, second period.}
\end{table}

\begin{table}
\centering\label{Tab:PortfolioObjectives2}
\makebox[\linewidth]{
\begin{tabu}{cccc}
\toprule
\   &   {$c_{2j}^1$} &   {$c_{2j}^2$} &   {$c_{2j}^3$}   \\
           \hline
\rowfont{\small}
\textit{$P_1$} 	&	60	&	145	&	42		\\
\textit{$P_2$}	&	43	&	54	&	365		\\
\textit{$P_3$}	&	75	&	150	&	57		\\
\bottomrule
\end{tabu}
}\caption{Objective function values, second period}
\end{table}

\noindent For the sake of the simplicity, we assume that the research salaries $w_{ij}$ are the same in the two periods as in Table 2.a. Similarly, we suppose that the $u_{jph}$ are the same for both periods as in Table 2.b.

After \textcolor{black}{generating the set of solutions in the Pareto Front}  with the $\epsilon-$constraint method we obtained the following twelve efficient portfolios of portfolios and outcome vectors.

\begin{enumerate}
    \item The project $P_1$ is selected and scheduled for period 1, and the project $P_2$ is selected for period 2, with variables $x_{11}=x_{22}=1$,  associated to the following possible combinations of researchers:
    \begin{itemize}[label={--}]
        \item researcher $e_1$ allocated to project $P_2$ and scheduled for period 2, researchers $e_2$ and $e_3$ allocated to project $P_1$ and scheduled for period 1 and researcher $e_4$ allocated to project $P_2$ and scheduled for period 2, with variables $y_{212},y_{121},y_{131},y_{242}$ equal to 1 and all others equal to 0;
    \end{itemize}
    \item The projects $P_1$ and $P_3$ are selected and scheduled for period 1 with variables $x_{11}=x_{13}=1$,  associated with the following possible combinations of researchers:
    \begin{itemize}[label={--}]
        \item researchers $e_2$, $e_3$ allocated to project $P_1$ and scheduled in period 1, researcher $e_4$ allocated to project $P_3$ and scheduled in period 1, with variables $y_{121},y_{131},y_{143}$ equal to 1 and all others equal to 0;
        \item researcher $e_2$ allocated to project $P_3$ and scheduled in period 1 and researchers $e_3$, $e_4$ allocated to project $P_1$ and scheduled in period 1, with variables $y_{123},y_{131},y_{141}$ equal to 1 and all others equal to 0;      
    \end{itemize}
    \item The project $P_1$ is selected and scheduled for period 1, and the project $P_3$ is scheduled for period 2 with variables $x_{11}=x_{23}=1$,  associated to the following possible combinations of researchers:
    \begin{itemize}[label={--}]
        \item researcher $e_2$ allocated to project $P_1$ and scheduled in period 1 and researchers $e_3$, $e_4$ allocated to project $P_3$ and scheduled in period 2, with variables $y_{121},y_{233},y_{243}$ equal to 1 and all others equal to 0;
        \item researchers $e_1$, $e_2$ allocated to project $P_3$ and scheduled in period 2 and researchers $e_3$ and $e_4$ allocated to project $P_1$ and scheduled in period 1, with variables $y_{213},y_{223},y_{141}$ equal to 1 and all others equal to 0;
    \end{itemize}
\item The  project $P_2$ is selected and scheduled for period 1, and  the project $P_3$ is scheduled for period 2
with variables $x_{12}=x_{23}=1$,  associated to the following possible combinations of researchers:
    \begin{itemize}[label={--}]
        \item researchers $e_1$, $e_2$ allocated to project $P_3$ and scheduled in period 2 and researchers $e_3$ and $e_4$ allocated to project $P_2$ and scheduled in period 1, with variables $y_{213},y_{223},y_{132},y_{142}$ equal to 1 and all others equal to 0;
        \item researchers $e_1$, $e_2$, $e_3$ allocated to project $P_3$ and scheduled in period 2 and researcher $e_4$ allocated to project $P_2$ and scheduled in period 1, with variables $y_{213},y_{223},y_{133},y_{142}$ equal to 1 and all others equal to 0;
        \item researchers $e_1$, $e_2$ allocated to project $P_3$ and scheduled in period 2 and researcher $e_4$ allocated to project $P_2$ and scheduled in period 1, with variables $y_{213},y_{233},y_{142}$ equal to 1 and all others equal to 0;
  \end{itemize} 
  \item The  project $P_2$ is selected and scheduled for period 1 and  the project $P_3$ is  scheduled for period 1 with variables $x_{12}=x_{13}=1$,  associated to the following possible combinations of researchers:
    \begin{itemize}[label={--}]
    \item researchers $e_1$ and $e_3$ allocated to project $2$ and scheduled in period 1, researchers $e_2$ and $e_4$ allocated to project $P_3$ and scheduled in period 1, with variables $y_{112},y_{132},y_{123},y_{143}$ equal to 1 and all others equal to 0;
    \item researchers $e_1$ and $e_2$ allocated to project $P_3$ and scheduled in period 1, researchers $e_3$ and $e_4$ allocated to project $P_2$ and scheduled in period 1, with variables $y_{113},y_{123},y_{132},y_{142}$ equal to 1 and all others equal to 0;
    \item researchers $e_1$, $e_2$  allocated to project $P_3$ and scheduled in period 1 and researcher $e_4$ allocated to project $P_2$ and scheduled in period 1, with variables $y_{113},y_{123},y_{142}$ equal to 1 and all others equal to 0;
    \item researchers $e_1$, $e_2$ and $e_3$ allocated to project $P_3$ and scheduled in period 1 and researcher $e_4$ allocated to project $P_2$ and scheduled in period 1, with variables $y_{113},y_{123},y_{133},y_{142}$ equal to 1 and all others equal to 0;
  \end{itemize}

\end{enumerate}

\subsection{Adding stochastic considerations}\label{subsec:stochastic_mo-portfolios}

\noindent Let $\Sigma = \{\sigma_{1},\ldots, \sigma_{\iota},\ldots,\sigma_{\tau}\}$ denote a set of states of nature, and $\pi_{\iota}$ denote the probability of state of nature $\sigma_{\iota}$; we have thus the set or probabilities, $\Pi = \{\pi_1,\ldots,\pi_{\iota},\ldots,\pi_{\tau}\}$. For each state of nature we have the performance levels, $c^{\ell}_{j}(\sigma_{\iota})$, with respect to each objective function, $z_\ell$, for $\ell=1,\ldots,k$. To evaluate the contributions of a portfolio, $\cal{P}$, we define the following expected value, for $\ell = 1,\ldots, k$ 

\begin{equation}
    {\displaystyle 
    {\cal{E}}({\cal{P}}) \;\;= \sum_{\{\iota\;:\;\sigma_{\iota} \in \Sigma\}}\sum_{j=1}^{n}c^{\ell}_{j}(\sigma_{\iota})x_j\pi_{\iota}.
    } 
\end{equation}

For each state of the nature, $\sigma_{\iota} \in \Sigma$, and each criterion, $g_p$, we also have a performance level for each element, $e_i \in E$, denoted by $g_p(e_i,\sigma_{\iota})$.  Inspired by the work of {Greco et al., 2010}, on the basis of the probability distribution, $\Pi$, we associated with each subset ${\cal{W}} \subseteq \Sigma$, the probability that one of the states is verified, i.e, ${\displaystyle \Pi_{\cal{W}}= \sum_{\{\iota \; : \; \sigma_{\iota} \in \Sigma \}} \Pi_{\iota}}$. Then, we can define a set $\Phi$ that includes all the potential probabilities that can be associated with any subset ${\cal{W}} \subseteq \Sigma$, as follows:

\begin{equation}
    {\displaystyle \Phi = \left\{ \phi \in [0,1] \; : \; \phi = \Pi_{\cal{W}}, \, {\cal{W}} \subseteq \Sigma \right\}.}
\end{equation}

Then, we can defined, for each criterion, $g_p$, a function $\alpha_p \; : \; E \times \Pi \rightarrow \Phi$, such that, for each element, $e_i \in E$, and each state, $\sigma_{\iota} \in \Sigma$, we have

\begin{equation}
    {\displaystyle 
    \alpha_p(e_i, \sigma_{\iota}) = \sum_{\{o\;:\; \sigma_o \in \Sigma,\; g_p(e_i,\sigma_o) \geqslant g_p{'}(e_i,\sigma_{\iota})\}}\pi_o.
    }
\end{equation}

The values $\alpha_p(e_i,\sigma_{\iota})$ represent the probability that element, $e_i$, obtains at least the performance level, $g_p{'}(e_i)$ for criterion $g_p$. Using the above function, we can define, for each criterion, $g_p$,  a function, $\rho \times \Phi \rightarrow \cal{X}$, as follows:

\begin{equation}
    {\displaystyle 
    \rho_p(e_i,\phi) = \max_{\{\iota \; : \; \sigma_{\iota} \Sigma, \; \alpha(e_i, \sigma_{\iota}) \geqslant \phi \}}\left\{g_p(e_i,\sigma_{\iota})\right\}.
    }
\end{equation}

The $\rho_p(e_i,\phi) = \cal{X}$ expresses that the evaluation, $g_p(e_i)$ is greater than or equal than $\cal{X}$ with probability at least equal to $\phi$. Thus, we can reformulate parameters, $v_{iph}$, as $v_{iph}(\phi)$  meaning that according to a probability value $\phi$, $v_{iph}(\phi)=1$ if $\rho_p(e_i, \phi) \geq l_s$, $v_{iph}(\phi)=0$, otherwise. Therefore, an element $e_i$ has a chance equal to $\phi$ to contribute to the satisfaction of the requirements expressed by the threshold levels $u_{jph}$.

Then, the \textit{sMO-PoP} can be formulated as follows for each $\phi$ value:
\begin{equation}\tag{\textit{sMO-PoP}}
    \begin{array}{lcl}
        {\displaystyle \max \; z_1(x_1,\ldots,x_j,\ldots,x_n) } \; {\displaystyle 
        {\cal{E}}({\cal{P}}) \;\;= \sum_{\{\iota\;:\;\sigma_{\iota} \in \Sigma\}}\sum_{j=1}^{n}c^{1}_{j}(\sigma_{\iota})x_j\pi_{\iota}.} \\
        \hspace{5.75cm} \vdots & & \\
         {\displaystyle \max \; z_\ell(x_1,\ldots,x_j,\ldots,x_n) } \; {\displaystyle 
        {\cal{E}}({\cal{P}}) \;\;= \sum_{\{\iota\;:\;\sigma_{\iota} \in \Sigma\}}\sum_{j=1}^{n}c^{\ell}_{j}(\sigma_{\iota})x_j\pi_{\iota}.} \\
        \hspace{5.75cm} \vdots & & \\
          {\displaystyle \max \; z_k(x_1,\ldots,x_j,\ldots,x_n) } \; {\displaystyle 
        {\cal{E}}({\cal{P}}) \;\;= \sum_{\{\iota\;:\;\sigma_{\iota} \in \Sigma\}}\sum_{j=1}^{n}c^{k}_{j}(\sigma_{\iota})x_j\pi_{\iota}.}\\
         \mbox{subject to:}  \\
            {\displaystyle \hspace{2.5cm} \sum_{j=1}^{n}y_{ij} \leqslant 1, \;\;\; i = 1, \ldots, m,}  \\
           {\displaystyle \hspace{2.5cm} y_{ij} \leqslant x_j, \;\;\; i = 1, \ldots, m, \; j = 1, \ldots, n,} \\
           {\displaystyle \hspace{2.5cm} \sum_{i=1}^{m}v_{iph(\phi)}y_{ij} \geqslant u_{jph}x_j, \;\;\; (j,p,h) \in D,} \\
           {\displaystyle \hspace{2.5cm} \sum_{i=1}^{m}\sum_{j=1}^{n}w_{ij}y_{ij} \leqslant W.} \\
    \end{array}
\end{equation}

\textcolor{black}{It is important to note that we have assumed the probabilities of an element reaching a specific threshold, as well as the thresholds themselves, are predetermined. However, we could also consider a scenario in which both the probabilities and the thresholds vary depending on the project allocations.}

\subsection{A third illustrative example}\label{subsec:ex_third_mo-portfolios}

\noindent Let us assume that we  have three different scenarios, $\Sigma=\sigma_1, \sigma_2, \sigma_3$ with  associated  probability $\pi_l={0.25,0.35,0.40}$, respectively. Let us also assume that the state of nature $\sigma_1$ corresponds to the deterministic scenario introduced in the first example. Instead, for the state of nature $\sigma_2$ and $\sigma_3$,  the evaluations $g_p(e_i,\sigma_\iota)$ are reported in Table \ref{tab:Criterion_Scenario2}. 
 
%  \begin{table}\caption{Objective Functions value   $c_{j}^l$ for scenarios $\sigma_2$ and $\sigma_3$}\label{tab:Criterion_Scenario}
% \centering
% \begin{subtable}{.50\textwidth}
% \centering
% \begin{tabular}{cccc}
% \hline
% \   &   {$F_1$} &   {$F_2$} &   {$F_3$}   \\
% \toprule
% \textit{$P_1$} 	&	25	&	22	&	30		\\
% \textit{$P_2$}	&	38	&	27	&	50		\\
% \textit{$P_3$}	&	47	&	27	&	38		\\
% \bottomrule
% \end{tabular}
% \caption{Scenario $\sigma_2$}
% \end{subtable}% <---- don't forget this %
% \begin{subtable}{.4\textwidth}
% \centering
% \begin{tabular}{cccc}
% \toprule
% \   &   {$F_1$} &   {$F_2$} &   {$F_3$}   \\
% \toprule
% \textit{$P_1$} 	&	50	&	24	&	47		\\
% \textit{$P_2$}	&	27	&	27	&	38		\\
% \textit{$P_3$}	&	26	&	41	&	23		\\
% \bottomrule
% \end{tabular}
% \caption{Scenario $\sigma_3$}
% \end{subtable}
% \end{table}

 \begin{table}\label{tab:Criterion_Scenario2}
\centering
\begin{subtable}{.50\textwidth}
\centering
\begin{tabular}{cccc}
\toprule
\    &   {$g_1$} &   {$g_2$} &   {$g_3$}   \\
       \hline
\textit{$e_1$} 	&	60	&	1	&	42		\\
\textit{$e_2$}	&	17	&	1	&	31		\\
\textit{$e_3$}	&	54	&	3	&	20		\\
\textit{$e_4$} 	&	42	&	2	&	30		\\
\bottomrule
\end{tabular}
\caption{Scenario $\sigma_2$}
\end{subtable}% <---- don't forget this %
\begin{subtable}{.4\textwidth}
\centering
\begin{tabular}{cccc}
\toprule
\    &   {$g_1$} &   {$g_2$} &   {$g_3$}   \\
       \hline
\textit{$e_1$} 	&	44	&	3	&	26		\\
\textit{$e_2$}	&	43	&	2	&	45		\\
\textit{$e_3$}	&	24	&	1	&	47		\\
\textit{$e_3$}	&	25	&	1	&	24		\\
\bottomrule
\end{tabular}
\caption{Scenario $\sigma_3$}
\end{subtable}
\caption{Performance levels $g_p(e_i, \sigma_2, \sigma_3)$ }\label{tab:Criterion_Scenario2}
\end{table}

Then, we can define the values $\rho_j(e_h, \phi)$ for each criterion $g_j \in G$. For the sake of the space, we report the $\rho(e_h, \phi)$ in Table 8.b only for the first criterion $g_1$. For example, the value $\rho_1(e_1, 0.25)=60$ means that the evaluations of element $e_1$ is at least $60$ with a probability of $25\%$.

% \begin{table}
% \centering\caption{Contribution $\rho_1(e_h, \phi)$ }\label{Tab:Portfolioprob}
% \makebox[\linewidth]{
% \begin{tabular}{lcccc}
% \hline
% \  $\phi$ &   {$e_1$} &   {$e_2$} &   {$e_3$}&   {$e_4$}   \\
%            \hline

% 0.25	&60 &43&54&42			\\
% 0.35	&60&43&54&42		\\
% 0.40	&	44&43&43&36		\\
% 0.60	&	44&24&43&36		\\
% 0.65	&	44&24&24&25		\\
% 0.75	&	44&17&24&25		\\
% 1.00	&	18&17&24&25		\\
% \hline
% \end{tabular}
% }
% \end{table}

On the basis of the above information the value of  $v_{hjs}$ can be calculated for each of the values assumed by $\phi$. As an example, in Table 9 we report those values obtained for $\phi=0.40$.

\begin{table}[htpb!]
\begin{minipage}{.3\linewidth}
    \centering 
    \begin{tabular}{ccccccccc}
    \toprule
    \multicolumn{1}{c}{} & \multicolumn{2}{c}{$g_1$} & \multicolumn{2}{c}{$g_2$} & \multicolumn{2}{c}{$g_3$}  \\
    \cmidrule(rl){2-3} \cmidrule(rl){4-5} \cmidrule(rl){6-7} 
    $v_{iph}$ & $l_{11}$ & $l_{12}$ & $l_{21}$ & $l_{22}$ & $l_{31}$ & $l_{32}$  \\
    \midrule
    $e_1$ & 1 & 1 & 1 & 1 & 1 & 1 \\
    $e_2$ & 1 & 1 & 1 & 0 & 1 & 1 \\
    $e_3$ & 1 & 1 & 0 & 0 & 1 & 1  \\   
    $e_4$ & 1 & 0 & 1 & 0 & 1 & 0  \\ 
    \bottomrule
    \end{tabular}
    \caption*{a. Parameters $v_{iph(\phi)}$, $\phi=0.4$.}
  \end{minipage}
  \hspace{4cm}
 \begin{minipage}{.3\linewidth}
    \centering
    \begin{tabular}{ccccc}
    \toprule
\  $\phi$ &   {$e_1$} &   {$e_2$} &   {$e_3$}&   {$e_4$}   \\
           \hline
0.25	&60 &43&54&42			\\
0.35	&60&43&54&42		\\
0.40	&	44&43&43&36		\\
0.60	&	44&24&43&36		\\
0.65	&	44&24&24&25		\\
0.75	&	44&17&24&25		\\
1.00	&	18&17&24&25		\\
    \bottomrule
    \end{tabular}
    \caption*{b. Contributions $\rho_1(e_h, \phi)$.}
  \end{minipage}%
  \caption{Performance levels of the elements and minimal qualitative levels required, second period.}

\end{table}

In this case, the value $v_{111}=1$ in the first row and the first column indicated that there is a $40\%$ probability that with respect to criterion $g_1$ the evaluation of element $e_1$ is at least $44$. 

We solved this model with the $\epsilon-$constraint method for each potential scenario. In the following, for the sake of the space, we only report the number of solutions obtained:

\begin{itemize}[label={--}]
\item For probability $\phi=0.25$ we obtained 41 possible combinations of non-dominated portfolios;
\item For probability $\phi=0.35$ we obtained 28 possible combinations of non-dominated portfolios;
\item For probability $\phi=0.4$ we obtained 12 possible combinations of non-dominated portfolios;
\item For probability $\phi=0.6$ we obtained 9 possible combinations of non-dominated portfolios;
\item For probability $\phi=0.65$ we obtained 4 possible combinations of non-dominated portfolios.
\end{itemize}

% \begin{itemize}
   
%     \item For probability $\phi=0.25$
%     \begin{itemize}
        
%     \item The portfolios of projects $P_1$ and $P_2$ is selected
%     with a total of 17 possible combinations of researchers, for example:
%     \begin{itemize}
%         \item researchers $e_1$ and $e_2$ allocated to project $P_1$ and researcher $e_3$ and researcher $e_4$ allocated to project $P_2$;
%          \item researchers $e_1$ and $e_2$ allocated to project $P_2$ and researcher $e_3$ and researcher $e_4$ allocated to project $P_3$;
%     \end{itemize}
%     \item The portfolios of projects $P_2$ and $P_3$ is selected
%     with the following possible combinations of researchers with a total of 18 possible combinations ofresearchers,s, as for example::
%     \begin{itemize}
%         \item researchers $e_1$ and $e_4$  allocated to project $P_2$ and researcher $e_2$ and $e_3$ allocated to project $P_3$;
%         \item researchers $e_1$ and $e_2$  allocated to project $P_2$ and researcher $e_4$ allocated to project $P_3$;
%     \end{itemize}
%     \item The portfolios of projects $P_1$ and $P_3$ is selected
%     with the following possible combinations of researchers with a total of 9 possible combinations of researchers, as for example::
%     \begin{itemize}
%         \item researchers $e_3$ and $e_4$  allocated to project $P_1$ and researcher $e_1$ and $e_2$ allocated to project $P_3$;
%         \item researchers $e_1$ and $e_4$  allocated to project $P_3$ and researcher $e_2$ $e_3$ allocated to project $P_3$;
%     \end{itemize}
% \end{itemize}
% \end{itemize}

Let us note that for probability $\phi =0.75$, we could not define any portfolios because the constraints on the requirements were not satisfied. Let us note that these obtained portfolios can be presented to the DM, who can conduct a scenario analysis to define his preferred portfolio. In the next Section, we will define an interactive procedure to handle such a process.

\section{An interactive methodology}\label{sec:methodology}
\noindent {We utilize the IMO-DRSA method \citep{gms_IMO}. In the IMO-DRSA framework, an interactive procedure is integrated within a multi-objective optimization process using the Dominance-Based Rough Set Approach (DRSA) \citep[see, e.g.,][]{greco2001rough, greco2010dominance}.}

{During the computation step, we solve a set of optimization problems to identify portfolios that will be presented to the DM. In the dialogue step, the DM can select a portfolio that she/he find completely satisfactory. If none meet this criterion, the DM is asked to evaluate which of the presented portfolios are relatively good. Following this evaluation, an algorithm based on DRSA generates ``if...,then..." rules to explain the DM's judgments in terms of minimal thresholds for the objective functions. For example, a decision rule might indicate that if a portfolio has an evaluation of at least 70 on the first objective to maximise and at least 50 on third objective to maximise, then it is considered relatively good. This implies that all portfolios meeting these thresholds are regarded as relatively good, and portfolios that do not meet the thresholds  are evaluated as not good.}

{The decision rules derived from this process are presented to the DM. At each iteration, the DM is tasked to select the rule she/he consider most representative of her/his preferences. The thresholds outlined in the selected rule are then added to the constraints of the problem, guiding the process in the subsequent iteration. This method effectively constrains the solution space to portfolios that satisfy all the thresholds established by the decision rules chosen by the DM, thereby supporting the selection of the final portfolio.}

{In the computational phase, we compute a set of non-dominated portfolios. For each portfolio obtained, we evaluate the values achieved by the objective functions, which can be discussed with the DM. For example, the DM may be interested in knowing how many projects are allocated to each portfolio.} 

{Various methods can be used to generate the non-dominated portfolios. One approach is to optimize a single objective in combinatorial optimization problems while ensuring that all the constraints of the original problem are maintained. After obtaining the portfolios, we can verify that they are non-dominated concerning the objectives that were not optimized and present these results to the DM.}

\textcolor{black}{We would like to emphasize that in this initial proposal for the new model, we chose to adopt an interactive approach. We believe this method is the most appropriate for addressing the problem, which inherently requires interaction with the decision-maker (DM). }

\section{llustrative case study with 8 researchers and 4 criteria}\label{sec:bigger_example}

\noindent{To illustrate the methodology, let us assume that we have a company that employs eight different researchers $E=\{e_1,\ldots, e_8\}$.
Following internal regulations, four criteria are used to evaluate each researcher, as follows:
\begin{enumerate}
    \item A score related to the years of experiences in managing R\&D projects of each researcher ($g_1$) ranging from 1 to 100;
    \item A score related to the expertise of the single researcher on a qualitative scale ranging from 1 to 6 ($g_2$); 
    \item A score related to the  number of previous projects completed  on time  ranging from 1 to 12 ($g_3$).
    \item A score related to the quality of previous innovations developed by the the researcher  ($g_4$), according to the following three levels qualitative scale:
    \begin{itemize}[label={--}]
        \item $0$, if no innovations have been developed;
        \item $1$, if innovations of low importance have been developed by the researcher;
        \item $2$, if innovations of great importance have been developed by the researcher;        
    \end{itemize}
\end{enumerate}
}

{
 In Table \ref{Tab:ProjectEvaluation}, we report the performance levels, $g_p(e_i)$,  on for each researcher across the four criteria. 
}

{
We identified five potential R\&D projects denoted as $P=\{P_1,P_2,P_3,P_4,P_5\}$,  each of them characterized by a minimum number of researchers required for every criterion. For example, 
 for criterion $g_1$, the score related to the years of experience in managing R\&D projects is  assessed based on two defined thresholds: a first threshold, $l_1(g_1)=20$, and a second threshold $l_1(g_1)=30$. Consequently, the requirements $v_{iph}$ can be specified. For instance, to select project $P_1$, we require a minimum of $3$ researchers with an evaluation for criterion $g_1$ of at least $20$ as well as at least $1$ researcher with an evaluation of at least  $30$. Therefore,  we can define the values $v_{111}=3$ and $v_{112}=1$. Similarly, for project  $P_2$  we need $2$ researchers with an evaluation for criterion $g_1$ of at least $20$ and at least $2$ researchers with an evaluation of at least $30$. Thus, we  define the values $v_{211}=2$ and $v_{212}=2$. Following this pattern, we can define the remaining values  $v_{iph}$. For brevity, these values are reported in the Appendix. 
}

\begin{table}[!h]
\centering\caption{Performance levels $g_p(e_i)$ for the illustrative case study.}\label{Tab:ProjectEvaluation}
\makebox[\linewidth]{
\begin{tabu}{ccccc}
\hline
\    &   {$g_1$} &   {$g_2$} &   {$g_3$} &   {$g_4$}   \\
       \hline
\rowfont{\small}

\textit{$e_1$} &	18&6&11&0\\
\textit{$e_2$}&	24&5&9&1\\
\textit{$e_3$}&	44&3&8&2\\
\textit{$e_4$} &	36&8&7&2\\
\textit{$e_5$} &	36&10&7&2\\
\textit{$e_6$} &	19&9&12&2\\
\textit{$e_7$} &12&5&4&2\\
\textit{$e_8$} &	26&4&8&2\\

\hline
\end{tabu}
}
\end{table}

{
Given the different required composition of each project, their contribution to each objective will be different. In Table \ref{Tab:CallsObjectives} we indicate those values, supposing that they are belong to a scale with a range [0,100]. The company evaluates the projects according to the following three strategic objectives: 
\begin{itemize}[label={--}]
    \item Scientific quality of the portfolio of projects ($c_j^1$);
    \item Adequacy to the organization strategy ($c_j^2$);
    \item Impact generated outside organisation ($c_j^3$).
\end{itemize}
}

{
In Table \ref{Tab:ProjectsCosts} we report the salary required by each researcher for each different project. The company has defined a maximum budget to be distributed among the projects equals to £300 millions.
}

\begin{table}[!h]
\centering\caption{Research salaries, $w_{ij}$}\label{Tab:ProjectsCosts}
\makebox[\linewidth]{
\begin{tabu}{ c c c c c c}
\hline
   &   {$w_{i1}$} &   {$w_{i2}$} &   {$w_{i3}$}&  {$w_{i4}$} &   {$w_{i5}$}  \\
       \hline
          
\rowfont{\small}
\textit{$e_1$} 	&	37	&	65	&	150	& 29 & 123	\\
\textit{$e_2$}	&	139	&	35	&	23 & 46&119		\\
\textit{$e_3$}	&	123	&	26	&	76	& 52&23	\\
\textit{$e_4$} 	&	119	&	28	&	14&25&76		\\
\textit{$e_5$} 	&	35	&	125	&	56&144&46		\\
\textit{$e_6$} 	&	24	&	35	&	63&186&52		\\
\textit{$e_7$} 	&	27	&	67	&	23&46&125		\\
\textit{$e_8$} 	&	98	&	86	&	57&136&35		\\
\hline
\end{tabu}
}
\end{table}
\begin{table}
\centering\caption{Objective function values}\label{Tab:CallsObjectives}
\makebox[\linewidth]{
\begin{tabu}{c c c c }
\hline
\   &   {$c_j^1$} &   {$c_j^2$} &   {$c_j^3$}   \\
           \hline
\rowfont{\small}
\textit{$P_1$} 	&	80	&	145	&	42		\\
\textit{$P_2$}	&	43	&	54	&	24		\\
\textit{$P_3$}	&	75	&	150	&	57		\\
\textit{$P_4$}	&	32	&	5	&	78		\\
\textit{$P_5$}	&	62	&	55	&	31		\\

\hline
\end{tabu}
}
\end{table}

\begin{table}[!h]
\centering\caption{     Minimal number of components $e_i$  assigned to each project, $u_{jph}$.}\label{Tab:MinimalNumberElementsGrande}
\makebox[\linewidth]{
\begin{tabu}{c c c c c c c c c c c}
\hline
\ \multirow{2}{*}{$u_{jph}$}   &  \multicolumn{2}{c} {$P_1$} & \multicolumn{2}{c}  {$P_2$} &   \multicolumn{2}{c} {$P_3$} & \multicolumn{2}{c}  {$P_4$} &   \multicolumn{2}{c} {$P_5$}  \\
       \cmidrule(rl){2-3} \cmidrule(rl){4-5} \cmidrule(rl){6-7} \cmidrule(rl){8-9} \cmidrule(rl){10-11}
      & $l_{11}$ & $l_{12}$ & $l_{21}$ & $l_{22}$ & $l_{31}$ & $l_{32}$ & $l_{41}$ & $l_{42}$ & $l_{51}$ & $l_{52}$   \\
      \hline
$g_1$ & 1 & 1 & 1 & 0 & 1 & 1& 1 & 1& 2 & 1 \\
$g_2$ & 1 & 1 & 1 & 2 & 2 & 1& 1 & 1& 1 & 0 \\
$g_3$ & 1 & 0 & 3 & 0 & 1 & 1& 2 & 1& 1 & 1 \\  $g_4$ & 2 & 0 & 3 & 1 & 1 & 0 & 0 & 0& 1 & 1\\ 
   
\hline
\end{tabu}
}
\end{table}

\subsection{An interactive methodology for the \textit{MO-PoP} model}

\noindent {To find the most preferred portfolio of portfolios, we introduce the following quantities:
\begin{itemize}[label={--}]
    \item $\chi=\sum\limits_{e_i \in E}\sum\limits_{P_j \in \mathcal{P}} y_{ij}$,  represents the total number of elements that  have been allocated to a portfolio $\mathcal{P}$;
      \item a set of threshold levels $L_q=\lambda_{11},\ldots,\lambda_{sl},\ldots, \lambda_{ml}$ for each objective $z_l$ for $ l=1,\ldots,k$;
     \item the number of projects in the portfolio $\mathcal{P}$ that have an evaluation at least equal to the level $\lambda_{sl}$ for objective $z_l$ and for level $\lambda_{sl} \in L_q$ as $\mathcal{F}_{s,\lambda}=|\sum\limits_{P_j \in \mathcal{P}}x_j: c_j^l\geqslant \lambda_{sl}|$.
\end{itemize}
}

{
To illustrate the methodology, we will apply it to the R\&D problem mentioned in the previous subsection. First, we will use the CPLEX optimization software to solve the single-objective problems by adopting each of the quantities $\mathcal{F}_{s\lambda}$ outlined earlier as the objective function, along with all the constraints defined for the \textit{MO-PoP} model.  The solutions are reported in Table \ref{Portfolios_first_iterationl}.  For instance, portfolio $\mathcal{P}_1$ is the optimal portfolio that maximizes the  number of projects $P_j$ for which the contribution $c_j^1$ is at least equal to the first qualitative threshold $\lambda_{11}$. 
}

{
Additionally, we could define an optimization problem using the quantities  $\chi$ and $\chi_j$ for each project $P_j \in \mathcal{P}$. However, in this particular case, it is more meaningful to minimize the number of researchers involved, which may render the problem more relevant.  We can engage in a conversation with the DM to determine if maximizing those objectives is the right approach. Subsequently, the portfolios can be presented to the DM. In our example, we simulated the role of the DM adopting an utility function. More in detail, 
we define a ``good portfolio" as one that meets a certain threshold based on a utility function that maximizes the sum of all contributions, denoted as $\mathcal{F}_l$. Therefore, all the portfolios but portfolio $\mathcal{P}_1$ are considered as ``good", as reported in the last column in Table \ref{Portfolios_first_iterationl}. 
}

\begin{table}
\centering 
 \caption{The set of non-dominated portfolios $\mathcal{P}$ obtained in the first iteration for the POP model}\label{Portfolios_first_iterationl}

\begin{tabular}{|c||c|c|c|c|c|c|c||c|}
%\hline
 %& \multicolumn{4}{c|}{\textbf{Criteria}}\\
%\hline
\hline

Portfolio&	{$\chi$}&		{$\mathcal{F}_{1,1}$}&	{$\mathcal{F}_{2,1}$}&	{$\mathcal{F}_{3,1}$}&{$\mathcal{F}_{1,2}$}&	{$\mathcal{F}_{2,2}$}&	{$\mathcal{F}_{3,2}$}&	
Evaluation\\
\hline

$\mathcal{P}_1$&8&4&1&2&3&1&2&*\\
$\mathcal{P}_2$&7&3&2&2&3&2&1&Good\\
$\mathcal{P}_3$&7&3&1&3&3&1&2&Good\\
$\mathcal{P}_4$&7&3&2&2&3&2&1&Good\\
$\mathcal{P}_5$&7&3&2&2&3&2&1&Good\\
$\mathcal{P}_6$&7&3&1&3&3&1&2&Good\\
\hline
\end{tabular}

\end{table}

{
At this point, we apply the DRSA procedure. We obtain a total of 5 rules, i.e.
\begin{itemize}[label={--}]
\item {\bf Rule 1.1:}  if $\chi_\leq 7$, then portfolio $\mathcal{P}$  is ``good",          
\item[~] (if there are less than seven researchers allocated to all the projects, then the portfolio is good).
	\item {\bf Rule 1.2:}  if $\mathcal{F}_{2,1}(\mathcal{P})\geqslant 2$, then portfolio $\mathcal{P}$  is ``good", \;\;\;\;   \;\;\;\;  \;\;\;\;         
\item[~] (if there are at least two projects with a contribution for objective $F_2$ at least equal to qualitative threshold $\lambda_{21}=30$, then the portfolio is good);
	\item {\bf Rule 1.3:}  if $\mathcal{F}_{2,2}(\mathcal{P})\geqslant 2$, then portfolio $\mathcal{P}$  is ``good", \;\;\;\;   \;\;\;\;  \;\;\;\;          
\item[~] (if there are at least two projects with a contribution for objective $F_2$ at least equal to the second threshold $\lambda_{22}=50$ at least equal to the second threshold, then the portfolio is good);
\item {\bf Rule 1.4:} if $\mathcal{F}_{3,1}(\mathcal{P})\geqslant 3$, then portfolio $\mathcal{P}$  is ``good", \;\;\;\;   \;\;\;\;  \;\;\;\;                
\item[~] (if there are at least three projects with a contribution for objective $F_3$ at least equal to the first threshold $\lambda_{31}=30$ , then the portfolio is good).
\end{itemize}
}

{
It is easy to verify that only Rule 1.1 is satisfied by all portfolios categorized as  ``good" in Table \ref{Portfolios_first_iterationl} , i.e, all the ``good" portfolios are characterised by having a total number of selected researchers smaller than 7.  
Specifically, all ``good" portfolios are characterized by having fewer than 7 selected researchers. At this stage, the rules can be presented to the DM, who can choose one or more of them to introduce as constraints in the optimization problem. 
}

{
Assuming that the DM chooses the rule that holds for all the portfolios, i.e. rule 1.1.,  we introduce the  constraint
$$\sum_{i=1}^{m}\sum_{j=1}^ny_{ij} \leqslant 7.$$ 
}

{
This modifies the single-objective models defined previously leading to the generation of a new set of portfolios is obtained. In this scenario, the optimal portfolios  are the same as portfolios $\mathcal{P}_2$ and $\mathcal{P}_4$ of the first iteration. Therefore, the DM is invited invited to choose between either of these portfolios, and the process concludes.
}

\subsection{An interactive methodology for the \textit{tMO-PoP} model}
\noindent {In order to find the most preferred portfolio of portfolios with temporal considerations, we add to the previous notation the following quantities:
\begin{itemize}[label={--}]
       \item the number of projects in portfolio $\mathcal{P}$ have an evaluation at least equal to the level $\lambda_{sl}$ for objective $z_l$ and for level $\lambda_{sl} \in L_q $ for period $t \in T$ as $\mathcal{F}_{t,s,\lambda}=|\sum\limits_{P_j \in \mathcal{P}}x_{tj}: c_{tj}^{l}\geqslant \lambda_{sl}|$.
\end{itemize}
}

\begin{table}
\centering \footnotesize \caption{The set of non-dominated portfolios $\mathcal{P}^T$ obtained in the first iteration for the POPT model}\label{Portfolios_first_iterationtempo}
\makebox[\linewidth]{
\begin{tabular}{|c||c|c|c|c|c|c|c||c||c|c|c|c|c|c|c|}
%\hline
 %& \multicolumn{4}{c|}{\textbf{Criteria}}\\
%\hline
\hline

Portfolio&	{$\chi$}&		{$\mathcal{F}_{1,1,1}$}&	{$\mathcal{F}_{2,1,1}$}&	{$\mathcal{F}_{3,1,1}$}&{$\mathcal{F}_{1,2,1}$}&	{$\mathcal{F}_{2,2,1}$}&	{$\mathcal{F}_{3,2,1}$}&	{$\mathcal{F}_{1,1,2}$}&	{$\mathcal{F}_{2,1,2}$}&	{$\mathcal{F}_{3,1,2}$}&{$\mathcal{F}_{1,2,2}$}&	{$\mathcal{F}_{2,2,2}$}&	{$\mathcal{F}_{3,2,2}$}	&
Evaluation\\
\hline

$\mathcal{P}^T_1$&6&3&3&3&3&3&1&0&0&0&0&0&0&*\\
$\mathcal{P}^T_2$&6&3&3&3&3&3&1&0&0&0&0&0&0&*\\
$\mathcal{P}^T_3$&6&3&3&3&3&3&1&0&0&0&0&0&0&*\\
$\mathcal{P}^T_4$&6&3&3&3&3&3&1&0&0&0&0&0&0&*\\
$\mathcal{P}^T_5$&6&3&3&3&3&3&1&0&0&0&0&0&0&*\\
$\mathcal{P}^T_6$&6&3&2&3&2&2&2&0&0&0&0&0&0&Good\\
$\mathcal{P}^T_7$&8&0&0&0&0&0&0&4&1&2&3&1&2&*\\
$\mathcal{P}^T_8$&7&1&1&1&1&1&0&2&2&1&2&2&0&Good\\
$\mathcal{P}^T_9$&7&0&0&0&0&0&0&3&1&3&3&1&2&*\\
$\mathcal{P}^T_{10}$&7&0&0&0&0&0&0&3&1&3&3&1&2&*\\
$\mathcal{P}^T_{11}$&7&1&1&1&1&1&0&2&2&1&2&2&0&Good\\
$\mathcal{P}^T_{12}$&7&0&0&0&0&0&0&3&1&3&3&1&2&*\\

\hline
\end{tabular}
}
\end{table}

{
First, we adopt each of the quantities $\mathcal{F}_{t,s,\lambda}$  introduced above  as objective function. We then use an optimisation software to solve the the single objective problems, obtaining the solutions reported in Table \ref{Portfolios_first_iterationtempo}.
By applying the DRSA approach, we identify 15 rules. None of these  is satisfied by all the ``good" portfolios; however two out of the three portfolios classified as ``good" satisfy all the rules. For the sake of brevity, we only present the  selected rule adopted in the procedure:
}

{
\begin{itemize}[label={--}]
	\item {\bf Rule $1.1^T$:} if $\mathcal{F}_{1,2,2}(\mathcal{P})\ge 1$ and $\mathcal{F}_{2,1,1}(\mathcal{P})\ge 2$ in the first period, then portfolio $\mathcal{P}$  is ``good", \;\;\;\;   \;\;\;\;  \;\;\;\;   
\item[] (if there is at least 1 project with a contribution for objective $F_1$ at least equal to qualitative threshold $l_{11}=20$ and if there are at least 2 projects with a contribution for objective $F_2$ at least equal to qualitative threshold $l_{11}=20$ in the first period, then the portfolio is good);
\end{itemize}
}

{
 We then recalculated  a new set of portfolios based on the newly introduced rules as shown in Table \ref{Portfolios_second_iterationtempo}. 
This time we identify nine rules, of which only two rules are satisfied by all the portfolios. In more detail:
}

\begin{table}
\centering \footnotesize \caption{The set of non-dominated portfolios $\mathcal{P}^{T'}$obtained in the second iteration for the POPT model}\label{Portfolios_second_iterationtempo}
\makebox[\linewidth]{
\begin{tabular}{|c||c|c|c|c|c|c|c||c||c|c|c|c|c|c|c|}
%\hline
 %& \multicolumn{4}{c|}{\textbf{Criteria}}\\
%\hline
\hline

Portfolio&	{$\chi$}&		{${F}_{1,1,1}$}&	{${F}_{2,1,1}$}&	{${F}_{3,1,1}$}&{${F}_{1,2,1}$}&	{${F}_{2,2,1}$}&	{${F}_{3,2,1}$}&	{${F}_{1,1,2}$}&	{${F}_{2,1,2}$}&	{${F}_{3,1,2}$}&{${F}_{1,2,2}$}&	{${F}_{2,2,2}$}&	{${F}_{3,2,2}$}	&
Evaluation\\
\hline

$\mathcal{P}_1^{T'}$&8&2&1&2&1&1&2&2&1&0&1&1&1&*\\
$\mathcal{P}_2^{T'}$&8&2&2&1&1&2&0&2&0&1&1&0&1&*\\
$\mathcal{P}_3^{T'}$&8&2&1&2&1&1&1&2&1&1&1&1&1&*\\
$\mathcal{P}_4^{T'}$&8&2&1&2&1&1&1&2&1&1&1&1&1&*\\
$\mathcal{P}_5^{T'}$&8&2&2&1&1&2&0&2&0&1&1&0&1&*\\
$\mathcal{P}_6^{T'}$&8&2&1&2&1&1&2&2&1&0&1&1&1&*\\
$\mathcal{P}_7^{T'}$&7&1&1&1&1&1&1&2&1&2&2&1&2&Good\\
$\mathcal{P}_8^{T'}$&7&1&1&1&1&1&0&2&2&1&2&2&2&Good\\
$\mathcal{P}_9^{T'}$&7&1&1&1&1&1&1&2&1&2&2&1&2&Good\\
$\mathcal{P}_{10}^{T'}$&7&1&1&1&1&1&1&2&1&2&2&1&2&Good\\
$\mathcal{P}_{11}^{T'}$&7&1&1&1&1&1&0&2&2&1&2&2&2&Good\\
$\mathcal{P}_{12}^{T'}$&6&1&1&1&1&1&0&2&0&2&2&0&2&Good\\

\hline
\end{tabular}
}
\end{table}

{
\begin{itemize}[label={--}]
\item  {\bf Rule $2.1^T$:} if $\chi_\leq 7$, then portfolio $\mathcal{P}$  is ``good",          
\item[~] (if there are less than 7 researchers allocated to all the projects, then the portfolio is good).
	\item {\bf Rule $2.2^T$:} if $\mathcal{F}_{1,2,2}(\mathcal{P})\ge 2$, then portfolio $\mathcal{P}$  is ``good", \;\;\;\;   \;\;\;\;  \;\;\;\;                
\item[~] (if there are at least 2 projects with a contribution for objective $F_2$ at least equal to qualitative threshold $l_{22}=30$ in the second period, then the portfolio is good);
\end{itemize}
}

{We select rule $2.1^T$, which is added as a constraint to our model and optimise the quantities $F_{\lambda,s,l}$. At this stage, only two potential portfolios are identified: $\mathcal{P}_{1}^{T'}$ and $\mathcal{P}_{12}^{T'}$. Thus, one of these two portfolios can be chose as the preferred solution. In a real case study the DM will select one of the two listed portfolios, at which point the procedure will conclude. Let us note that the solutions proposed to the DM can be calculated using various methods. One such approach could be the $\epsilon$-constraint method suggested by \citep{mesquita2023new},  which allows us to identify an even larger number of solutions. Depending on the specific problem and the preferences of the DM, we could propose a greater or smaller number of solutions. In this example, we have chosen to optimize each objective function individually to ensure a diverse set of solutions.  Furthermore, the representation problem of the Pareto front involves determining how many portfolios to present to the decision-maker (DM) and evaluating how well those portfolios cover and represent the Pareto front. Among the latest methods we can cite the GRid Point Based Algorithms \citep{mesquita2023new} and   the exact method proposed by \citep{fotedar2025method}. We plan to investigate this aspect in future research. 
}

\subsection{An interactive methodology for the \textit{sMO-PoP model}}
 \noindent {Let us assume that we  have three different scenarios, with  associated  probability $\Pi={0.25,0.35,0.40}$. For the sake of the space, we report in the Appendix the evaluations for each of these scenarios. Then, we can define the $\rho(e_i, \phi)$ for each criterion $g_p$.  Then, the  $v_{iph}$ can be calculated for each of the values assumed by $\phi$. We selected a certain number of specific $\phi$ values and for each of them we evaluate the following quantities:
 \begin{itemize}[label={--}]
 \item $F_{\lambda,s}^{\phi}$ represent the number of projects in portfolio $\mathcal{P}$ have an evaluation at least equal to the level $\lambda_{sl}$ for objective $z_l$ and for level $\lambda_{sl} \in L_q$  with a probability at least equals to $\phi$ $F_{s,\lambda}^{\phi}=|\sum\limits_{P_j \in \mathcal{P}}x_j: c_j^l\geq \lambda_{sl}|$ with a probability at least equals to $\phi$.
\end{itemize}
}

{
As an example, let us describe a situation in which two values of $\phi$ are selected, e.g. 0.4 and 0.65. We start optimising the monobjective problems maximising each of the quantities $F_{s,\lambda}^{\phi}$ for each  objective $z_l$ and for each level $\lambda_{sl} \in L_q$. The obtained portfolios $\mathcal{P}^\phi$ are reported in Table \ref{Portfolios_first_iterationstoc}. 
}

\begin{table}
\centering \footnotesize \caption{The set of non-dominated obtained portfolios $\mathcal{P}^S$ in the first iteration.}\label{Portfolios_first_iterationstoc}
\makebox[\linewidth]{
\begin{tabular}{|c||c|c|c|c|c|c|c|c|c|c|c|c|c||c|}
%\hline
 %& \multicolumn{4}{c|}{\textbf{Criteria}}\\
%\hline
\hline

Portfolio&	{$\chi$}&		{${F}_{1,1}$}&	{${F}_{2,1}^{0.4}$}&	{${F}_{3,1}^{0.4}$}&{${F}_{1,2}^{0.4}$}&	{${F}_{2,2}^{0.4}$}&	{${F}_{3,2}^{0.4}$}&	{${F}_{1,1}^{0.65}$}&	{${F}_{2,1}^{0.65}$}&	{${F}_{3,1,s_3}^{0.65}$}&{${F}_{1,2}^{0.65}$}&	{${F}_{2,2}^{0.65}$}&	{${F}_{3,2}^{0.65}$}&Evaluation	\\
\hline

$\mathcal{P}_1^S$&8&4&4&4&4&3&1&3&2&2&2&2&1&*\\
$\mathcal{P}_2^S$&8&4&4&4&4&3&1&3&2&2&2&2&1&*\\
$\mathcal{P}_3^S$&8&4&4&4&4&3&1&3&2&2&2&2&1&*\\
$\mathcal{P}_4^S$&8&4&4&4&4&3&1&3&2&2&2&2&1&*\\
$\mathcal{P}_5^S$&8&4&4&4&4&3&1&3&2&2&2&2&1&*\\
$\mathcal{P}_6^S$&6&3&3&3&3&2&2&3&2&2&2&2&1&Good\\
$\mathcal{P}_7^S$&6&3&3&3&3&2&2&3&2&2&2&2&1&Good\\
$\mathcal{P}_8^S$&6&3&3&3&3&3&1&3&3&1&3&3& &Good\\
$\mathcal{P}_9^S$&6&3&3&3&3&2&2&3&2&2&2&2&1&Good\\
$\mathcal{P}_{10}^S$&6&3&3&3&3&3&1&3&3&1&3&3&0&Good\\
$\mathcal{P}_{11}^S$&6&3&3&3&3&3&1&3&3&1&3&3&0&Good\\
$\mathcal{P}_{12}^S$&6&3&3&3&3&2&2&3&2&2&2&2&1&Good\\

\hline
\end{tabular}
}
\end{table}
{
Applying the DRSA approach we obtained 5 rules. Only two of these rules obtained were satisfied by all the ``good" portfolios,i.e.
\begin{itemize}[label={--}]
	\item {\bf Rule $1.1^S$:} if $\mathcal{F}_{3,1}^{0.65}(\textbf{x})\ge 2$, then portfolio $\mathcal{P}$  is ``good", \;\;\;\;   \;\;\;\;  \;\;\;\;                
\item[~] (if with a probability of at least 65\% there are at least 2 researchers with a contribution for objective $F_3$ at least equal to qualitative threshold $l_{11}=20$, then the portfolio is good);
\item  {\bf Rule $1.^S$:} if $\mathcal{F}_{1,2}^{0.65}(\textbf{x})\ge 2$, then portfolio $\mathcal{P}$  is ``good", \;\;\;\;   \;\;\;\;  \;\;\;\;                
\item[] (if with a probability of at least 65\% there are at least 2 researchers with a contribution for objective $F_1$ at least equal to qualitative threshold $l_{11}=30$, then the portfolio is good);
\end{itemize}
}

{
We suppose that the DM considered as most representative of his preferences  and we calculated again a new set of portfolios according to the newly introduced rule; those are reported in Table \ref{Portfolios_second_iterationstoc}. 
}

\begin{table}
\centering \footnotesize \caption{The set of non-dominated obtained portfolios $\mathcal{P}^{S'}$ in the second iteration.}\label{Portfolios_second_iterationstoc}
\makebox[\linewidth]{
\begin{tabular}{|c||c|c|c|c|c|c|c|c|c|c|c|c|c||c|}
%\hline
 %& \multicolumn{4}{c|}{\textbf{Criteria}}\\
%\hline
\hline

Portfolio&	{$\chi$}&		{${F}_{1,1}^{0.4}$}&	{${F}_{2,1}^{0.4}$}&	{${F}_{3,1}^{0.4}$}&{${F}_{1,2}^{0.4}$}&	{${F}_{2,2}^{0.4}$}&	{${F}_{3,2}^{0.4}$}&	{${F}_{1,1}^{0.65}$}&	{${F}_{2,1}^{0.65}$}&	{${F}_{3,1,s_3}^{0.65}$}&{${F}_{1,2}^{0.65}$}&	{${F}_{2,2}^{0.65}$}&	{${F}_{3,2}^{0.65}$}&Evaluation	\\
\hline

$\mathcal{P}_1^{S'}$&8&4&4&4&4&3&1&3&2&2&2&2&1&Good\\
$\mathcal{P}_2^{S'}$&8&4&4&4&4&3&1&3&2&2&2&2&1&Good\\
$\mathcal{P}_3^{S'}$&8&4&4&4&4&3&1&3&2&2&2&2&1&Good\\
$\mathcal{P}_4^{S'}$&8&4&4&4&4&3&1&3&2&2&2&2&1&Good\\
$\mathcal{P}_5^{S'}$&8&4&4&4&4&3&1&3&2&2&2&2&1&Good\\
$\mathcal{P}_6^{S'}$&6&3&3&3&3&2&2&3&2&2&2&2&1&*\\
$\mathcal{P}_7^{S'}$&6&3&3&3&3&2&2&3&2&2&2&2&1&*\\
$\mathcal{P}_8^{S'}$&9&4&4&4&4&3&1&3&2&2&2&2&1&*\\
$\mathcal{P}_9^{S'}$&6&3&3&3&3&2&2&3&2&2&2&2&1&*\\
$\mathcal{P}_{10}^{S'}$&9&4&4&4&4&3&1&3&2&2&2&2&1&Good\\
$\mathcal{P}_{11}^{S'}$&9&4&4&4&4&3&1&3&2&2&2&2&1&Good\\
$\mathcal{P}_{12}^{S'}$&6&3&3&3&3&2&2&3&2&2&2&2&1&*\\

\hline
\end{tabular}
}
\end{table}

{
By applying the DRSA approach, we identified five rules that are satisfied by all the good portfolios. For brevity, we are only reporting some of the rules that the decision-maker considered to be the most representative of their preferences.
\begin{itemize}[label={--}]
	\item {\bf Rule $2.1^S$:} if $\mathcal{F}_{1,2}^{0.40}(\textbf{x})\ge 4$, then portfolio $\mathcal{P}$  is ``good", \;\;\;\;   \;\;\;\;  \;\;\;\;                
\item[] (if with a probability of at least 40\% there are at least 4 researchers with a contribution for objective $F_2$ at least equal to qualitative threshold $l_{11}=20$, then the portfolio is good);
\end{itemize}
}

{
After recalculating according to the newly introduced rule, we found that no new portfolios were identified. Only two portfolios were obtained, corresponding to existing portfolios $\mathcal{P}_1^{S'}$ and $\mathcal{P}_6^{S'}$. As a result, the interaction can cease by asking the DM which portfolios he prefers.
}

\section{Conclusions}\label{sec:conclusion}
\noindent {In this paper, we introduced a new approach to portfolio management that facilitates the selection of a portfolio of projects based on the selection and allocation of elements to each project. Our model considers both temporal and stochastic factors and can be easily customized for a variety of applications. We also provided an illustrative example to demonstrate how a multi-objective optimization technique can be applied to this model.}

{We believe that our models represent a significant addition to the portfolio decision analysis theory, enriching this growing area of research with a newly developed framework. From a practical standpoint, these models can serve as valuable tools for DMS, enabling them to apply the models to various real-life case studies. They provide a way to describe complex situations in an intuitive manner. The interactive process allows DMs to collaboratively construct their preferred portfolio until they are completely satisfied with the results.}

{In future research, we envision implementing additional techniques to address multi-objective problems, particularly heuristic and metaheuristic approaches, especially in the context of larger dimensional problems. Specifically, we aim to develop an interactive procedure that integrates the use of the deck of cards method \citep{garcia2024deck} for collecting the preferences of decision-makers (DMs), which could help reduce the cognitive burden on them \cite{barbati2024deck}. At the same time, we will focus on improving the representation of the Pareto front by adopting exact algorithms for its definition \cite{mesquita2023new}. Additionally, the model could be expanded to include further considerations, such as synergies among projects or interdependencies among various elements \citep{liesio2021portfolio}. Furthermore, it would be valuable to apply the models to real-world decision-making scenarios to evaluate their performance under more challenging conditions.}

\section*{Acknowledgements}
\addcontentsline{toc}{section}{\numberline{}Acknowledgements}
\noindent {Jos\'e Rui Figueira has been supported by Portuguese national funds through the FCT - Foundation for Science and Technology, I.P., grant number UIDB/00097/2020.}

\newpage

\bibliographystyle{model2-names}
\bibliography{biblio}
%\addcontentsline{toc}

\end{document}